\documentclass{tran-l}

\usepackage{amssymb}

\usepackage{graphicx}
\usepackage{mathtools}

\newtheorem{theorem}{Theorem}[section]

\theoremstyle{definition}

\theoremstyle{remark}

\numberwithin{equation}{section}

\begin{document}

\title{Rank One Approximation as a Strategy for Wordle}


\author{Michael Bonthron}
\address{DePaul University -- Chicago, Illinois}
\curraddr{}
\email{mbonthro@depaul.edu}
\thanks{I would like to thank Dr. Ilie Ugarcovici (DePaul Unviersity) for advice and encouragement on the project, Amrita Bhattacharyya for editorial insight, and DePaul University for research support.}



\date{March 20, 2022}

\dedicatory{}
\begin{abstract}
This paper presents a mathematical method of playing the puzzle game Wordle. 
In Wordle, the player has six tries to guess a secret word. 
After each guess the player is told how their guess compares to the secret word. With the available information the player makes their next guess. 
This paper proposes combining a rank one approximation and latent semantic indexing to a matrix representing the list of all possible solutions.
Rank one approximation finds the dominant eigenvector of a matrix of words, and latent semantic indexing reveals which word is closest to the dominant eigenvector.
The word whose column vector is closest to the dominant eigenvector is chosen as the next guess. With this method the most representative word of the set of all possible solutions is selected. This paper describes how a word can be converted to a vector and the theory behind a rank one approximation and latent semantic indexing. 
This paper presents results demonstrating that with an initial guess of ``SLATE'' the method solves the puzzle in $4.04$ guesses on average, with a success rate of $98.7\%$.
\end{abstract}

\maketitle


\section{INTRODUCTION}
In recent months, Wordle, an online daily puzzle game, has taken over group messages and Facebook posts with green, yellow, and gray squares proclaiming one's victory over the daily challenge. 
The goal of Wordle is to guess a secret five letter word.
After each guess, the commonalities between the secret word and the guessed word are revealed.
With the available information, the player narrows down the list of possible solutions. 
To win the game, the player must guess the secret word within six tries, although to win bragging rights, the player must guess the secret word in fewer tries than their friends. 
Questions such as ``what is the best starting word?'' and ``what should the next guess be?'' motivate an analytical approach to the game. 
This paper responds to both questions by utilizing a rank one approximation in combination with latent semantic indexing. 
A rank one approximation will be used to find an idealized best guess from a list of words.
Next, latent semantic indexing will be applied to find the English word closest to the idealized best guess.
The closest word will be used as the next guess in the game to continually guess a word representative of the set.
The method, theory, and application of a rank one approximation of a matrix and latent semantic indexing will be discussed. 
Results will demonstrate that this method, when using an initial guess of ``SLATE'', solves the puzzle in $4.04$ guesses on average with a win rate of $98.7\%$.

\section{Game play}
The best way to understand Wordle is to play the game available online for oneself  \cite{wordlesite}, but a short description of the game is included for those not familiar. 

Wordle is played against a computer and begins with the computer selecting a five letter secret word from the English dictionary. 
The player guesses an initial word.
The computer responds by coloring each letter of the guessed word either green, yellow, or gray based on how each letter occurs in the secret word. 
Letters colored green indicate the letter occurs in that position in the secret word.
A letter colored yellow indicates the letter appears in the secret word, but not in the same position as it was guessed.
A letter colored gray indicates the letter does not appear in the secret word at all.
For example, if the secret word was ``MATHS'' and the player guessed ``MARCH'', ``M'' and ``A'' would be green, ``H'' would be colored yellow, and ``R'' and ``C'' would be gray.

Double letters require some nuance. 
If the guessed word includes double letters but the secret word only contains the letter once, the computer will only credit one occurrence of the letter with either a green or yellow color.
For example, if the secret word was “LAKES” and the player guessed “LLAMA”, the first “L” would be colored green, but the second “L” would be colored gray, representing that the solution only has one “L”. 
Conversely, if the secret word has a double letter but the guessed word only has one occurrence of the letter, the letter will be colored either green or yellow. 
The computer will not color the letter half-green-half-yellow to inform the player that there is a second occurrence of the letter.

After the computer's response, the player uses the information to choose their next guess. 
The process is repeated until either the player successfully guesses the secret word or runs out of guesses. 
To win the game, the player must correctly guess the secret word within six guesses, but winning in few moves is more prestigious. 

The game can be played in two modes: ``hard mode'' or ``easy mode.'' In hard mode, all information available \textbf{must} be used in each guess, while in easy mode this constraint is not added. 
In easy mode, players can take advantage of this lack of constraints, particularly when four out of the five letters of the secret word have been solved, and there are multiple letters remaining that could satisfy the final spot. 
This paper will only consider playing the game in hard mode.  

From the player’s perspective there are two aspects to the game: interpreting the response from the computer to create a set of possible solutions and selecting a word from the list of possible solutions as the next guess. 
The first aspect becomes an easy task with the implementation of computer code to sort through the dictionary and test each word against the available conditions.
The second question opens the door for mathematical analysis and discussion.
This paper will focus on the question of selecting the next guess from a list of possible solutions.

\section{Mathematical models}
As the game grew to a viral sensation, a number of mathematical methods were proposed as strategies for Wordle.
Sources have used frequency analysis of letters to suggest words starting with ``S'' make good first guesses \cite{aops} \cite{tomsguide}.
The popular video creator 3Blue1Brown used information theory to minimize entropy to play Wordle in easy mode \cite{3blue1brown}. 
This paper proposes a combination of a rank one approximation and latent semantic indexing to select the  most representative word from the set of possible solutions. 
The two concepts will be combined as a strategy for Wordle.
First, a single vector representing a matrix will be produce using a rank one approximation.
Second, latent semantic indexing will choose the word from a list which is closest to the representative vector. 

\subsection{Converting a word to a column vector}
Before applying any rank one approximations or latent semantic indexing, a list of words must be converted to matrix form. 
This section will describe two methods of converting a five letter word into a vector. 
The first method will consider how many times each letter occurs in the word while the second method will consider both the occurrence and position of each letter.

Consider a $26 \times 1$ column vector whose $i^{th}$ row represents the $i^{th}$ letter of the alphabet.
A five letter word can be converted to a vector with the value in the $i^{th}$ row indicating the number of times the $i^{th}$ letter of the alphabet occurs in the word. 
\begin{equation}
    \textit{word} := 
    \begin{bmatrix}
        a_{1,1}  \\
        a_{2,1}  \\
        a_{3,1}  \\

        \vdots  \\
        a_{25,1} \\
        a_{26,1} \\
    \end{bmatrix}
        \begin{matrix}
        A \\
        B \\
        C \\
        \vdots \\
        Y \\ 
        Z \\
    \end{matrix}
\end{equation}
For example, the word ``ABBOT'' can be converted to a $26 \times 1$ vector as:
\begin{equation}
    \textit{``ABBOT''}:=\begin{bmatrix}
        1 \\
        2 \\
        0 \\
        0 \\
        0 \\
        \vdots \\
        1 \\
        \vdots \\
        1 \\
        \vdots \\
        0 \\
        0 \\
    \end{bmatrix}
    \begin{matrix}
        A \\
        B \\
        C \\
        D \\
        E \\
        \vdots \\
        O \\
        \vdots \\
        T \\
        \vdots \\
        Y \\ 
        Z \\
    \end{matrix}
\end{equation}

With this method, words are not uniquely represented since any words which are anagrams will produce the same vector.

To consider both the frequency and position of each letter in a word, a more elaborate scheme is needed.
A five letter word can be described using a $26 \times 5$ scheme where the $j^{th}$ column represents the $j^{th}$ letter in the word and the $i^{th}$ row represents the $i^{th}$ letter in the alphabet. Any word can uniquely be described by the $26 \times 5$ structure as follows:

\begin{equation}
    \textit{word} := 
    \begin{bmatrix}
        a_{1,1} & a_{1,2} &  a_{1,3} & a_{1,4} & a_{1,5} \\
        a_{2,1} & a_{2,2} &  a_{2,3} & a_{2,4} & a_{2,5} \\
        a_{3,1} & a_{3,2} &  a_{3,3} & a_{3,4} & a_{3,5} \\

        \vdots & \vdots & \vdots  & \vdots \\
        a_{25,1} & a_{25,2} &  a_{25,3} & a_{25,4} & a_{25,5} \\
        a_{26,1} & a_{26,2} &  a_{26,3} & a_{26,4} & a_{26,5} \\
    \end{bmatrix}
        \begin{matrix}
        A \\
        B \\
        C \\
        \vdots \\
        Y \\ 
        Z \\
    \end{matrix}
\end{equation}

A value of $1$ in the location $a_{i,j}$ represents the $i^{th}$ letter of the alphabet occurring in the $j^{th}$ position of the word. All other indices will be $0$.
This structure can be manipulated into a $130 \times 1$ column vector by stacking the columns atop one another. The first $26$ rows of the vector correspond to the first letter of the word. The $27^{th}$ through $53^{rd}$ rows correspond to the second letter of the word and so on. Manipulating the $26 \times 5$ structure results in the following:

\begin{equation}
 \textit{word} :=
\begin{bmatrix}
a_{1,1}\\
a_{1,2}\\
\vdots \\
a_{1,26}\\
a_{2,1} \\
\vdots \\
a_{2,26} \\
a_{3,1} \\
\vdots \\
a_{5,1}\\
a_{5,2}\\
\vdots \\
a_{5,26}\\
\end{bmatrix}
\end{equation}

For example, to convert the word ``ABBOT'' to a column vector, the first row will be $1$ since there is an ``A'' in the first position. The $28^{th}$ and $52^{th}$ rows will be $1$ since the second and third letters are ``B.'' Determining the rows representing ``O'' and ``T'' is challenging to do mentally, but are represented by $1$'s in the $93^{rd}$ and $124^{th}$ rows. All other rows in the vector are $0$. This method represents any five letter word as a unique  $130 \times 1$ column vector.

If a single word is represented by a column vector, a list of $n$ words can be represented by an $m \times n$ matrix $A$ where each column in $A$ represents an individual word.
\begin{equation}
    A := \begin{bmatrix}
       W & W & & W\\
       O & O & & O \\
       R & R & \hdots & R\\
       D & D & & D\\
       1 & 2 & & n\\
    \end{bmatrix}
\end{equation}

After forming the matrix $A$, a rank one approximation and latent semantic indexing can be applied to select the word which is most representative of the set. 

\subsection{Rank one approximation}
A rank one approximation of a matrix approximates the matrix using only two vectors and a scaling factor. A rank one approximation begins by performing single value decomposition to a matrix $A$. Single value decomposition factors $A$, an $m\times n$ matrix, into the following form:
\begin{equation}
    A=USV^T
\end{equation}
where:

    $A := m \times n$ matrix
    
    $U := m \times m$ matrix whose columns are the left singular vectors $u_i$

    $S := m \times n$ matrix whose diagonal entries are the singular values $s_i$
    
    $V := n \times n$ matrix whose columns are the right singular vectors $v_i$
    
And where:

    $\{ u_1, u_2 \dots u_m \}$ is an orthonormal set of eigenvectors of $AA^T$
    
    $\{ s_1, s_2 \dots s_m \}$ is the set of eigenvalues of $AA^T$
    
    $\{ v_1, v_2 \dots v_n \}$ is an orthonormal set of eigenvectors of $A^TA$
    
The eigenvectors for both $U$ and $V$ are arranged in descending order based on their corresponding eigenvalues. If $m<n$ this is to say:

\begin{equation}
    A =\begin{bmatrix}
    u_{11}       & \dots & u_{1m} \\
    \vdots       & \ddots      & \vdots \\
    u_{m1}       & \dots & u_{mm}
\end{bmatrix}
\begin{bmatrix}
    s_{1}    & 0      & \dots   & 0      & \dots & 0\\
      0      & s_{2}  & \dots   & 0      & \dots & 0 \\
    \vdots   & \vdots & \ddots  & \vdots & \dots &\vdots \\
    0        & 0      & \dots   & s_{m}  & \dots & 0 \\
    \end{bmatrix}
\begin{bmatrix}
    v_{11}       & \dots & v_{1n} \\
    v_{21}       & \dots & v_{2n} \\
    \vdots       & \ddots      & \vdots \\
    v_{n1}       & \dots & v_{nn}
\end{bmatrix}^T
\end{equation}

Where:
\begin{equation}
    s_{1}>s_{2}>...>s_{m}
\end{equation}

and both
\begin{equation}
    u_i = \begin{bmatrix}
        u_{1i} \\
        u_{2i} \\
        \vdots \\
        u_{mi}
    \end{bmatrix}
\end{equation} 

and
\begin{equation}
    v_i^T=\begin{bmatrix}
        v_{1i} & v_{2i} & ... & v_{ni}
    \end{bmatrix}
\end{equation} 
correspond to the singular value $s_i$.

A rank one approximation of $A$ uses the singular value $s_1$ and the corresponding vectors $u_1$ and $v_1$ to approximate the matrix $A$ \cite{sauer}. The best approximation of $A$ using only two vectors and a singular value is:
\begin{equation}
    A\approx u_{1}s_{1}v_{1}^T
\end{equation}

If $u_1$ is the column vector that best represents $A$, a vector of interest is the column from $A$ which is closest to $u_1$. Latent semantic indexing is used to find the column from the original matrix that is closest to $u_1$.

\subsection{Latent semantic indexing}
Latent semantic indexing is used in search engines to find results which closely match a searched term.
Given a query vector, latent semantic indexing is a systematic way of sorting a set of vectors based on their similarities to the query \cite{LSI}.
For our application, the query vector will be $u_1$, the dominant eigenvector of $AA^T$.
We are interested in finding which column of $A$ is closest to $u_1$.
Cosine similarity is used to calculate the angle between each column of $A$ and $u_1$.
The angle $\theta_n$ between $u_1$ and the $n^{th}$ column of $A$, $a_n$, is determined by: 
\begin{equation}
    \theta_n = cos^{-1}(\frac{u_1^Ta_n}{||u_1|| ||a_n||})
\end{equation}

A value of $\theta$ close to $0$ indicates that the column is a close match to the query. 
Values of $\theta$ close to $\frac{\pi}{2}$ represent columns that are not a close match to the query.
The column of $A$ which makes the smallest angle to $u_1$ will be described as the most representative column from the set since it is nearest in direction to the dominant eigenvector of the matrix representing the set.

If multiple words make the same angle with the dominant eigenvector, this method cannot distinguish which word is a better guess. 
In this case, one of the words with the smallest value for $\theta$ will be randomly selected as the next guess.

\subsection{Application to Wordle}
The following steps are used to apply these concepts to Wordle:
\begin{enumerate}
    \item Create a list of possible solutions using the available information about the secret word.
    \item Using either of the described methods, convert the list of possible solutions into the matrix $A$ with each column of $A$ representing one possible solution. 
    \item Determine $u_1$, the dominant eigenvector of $AA^T$.
    \item Determine $\theta_n$, the angle between $u_1$ and the $n^{th}$ column of $A$.
    \item Find the word corresponding to the smallest value of $\theta_n$.
    \item Use the word from step 5 as the next guess since it is most representative of the set of possible solutions. Repeat the process until the secret word is correctly solved.
    \end{enumerate}

\subsection{Example}
Consider the following set of words:

\begin{center}
    ``CLUMP'', ``CLAMP'', ``RUNNY'', ``UNDER'', ``CAMPS'', ``CRUNK''
\end{center}

If you had to choose one word to represent the entire set, which word would you select?
Justification could argue for multiple words, perhaps ``CLUMP'' or ``CRUNK'' or ``CLAMP.''
Even with this small list of six words, it can still be unclear which word to choose.
Let's see how the previously described method can be used to make an informed selection.

First, we need to choose which method of converting a word into a vector will be used.
Using the method which considers the letter's location for each word, a $130 \times 1$ vector to represent each word is created.
The subsequent matrix, $A$, which represents the entire set is:
\begin{equation}
    A = \begin{bmatrix}
        0 & 0 & 0 & 0 & 0 & 0 \\
        0 & 0 & 0 & 0 & 0 & 0 \\
        1 & 1 & 0 & 0 & 1 & 1 \\
       \vdots & \vdots & \vdots & \vdots & \vdots & \vdots \\
        0 & 0 & 0 & 0 & 0 & 0 \\
    \end{bmatrix}
\end{equation}

Whose dominant eigenvector, $u_1$, is:
\begin{equation}
    u_1=\begin{bmatrix}
        0 \\
        0 \\
        0.592 \\
        \vdots \\
        0.0227 \\
        0 \\
    \end{bmatrix}
\end{equation}

The vector $u_1$ is sparse, containing only $17$ non-zero rows.
Ideally, $u_1$ could be used as the next guess, but with this method there is no interpretation for a word with $0.592$ C's as the first letter.
Rather, we find the vector closest to $u_1$ using latent semantic indexing (table 1).

If instead, we converted each word to a vector with the method that did not consider letter location, $A$ would be a  $26 \times 5$ matrix.
All the same, the dominant eigenvector is found and latent semantic indexing is used to find the angle for each word (table 1).

\begin{table}[ht]
\caption{Comparing each word's angle to the dominant eigenvector of their corresponding matrix for both the methods of converting a word into a column vector.}\label{resulttab}
\renewcommand\arraystretch{1.5}
\noindent\[
\begin{array}{|c|c|c|c|}
\hline
\multicolumn{2}{|c|}{\textrm{Method not considering letter location}} & \multicolumn{2}{|c|}{\textrm{Method considering letter location}}\\
\hline
1 & \textrm{``CRUNK'' at } 36^{\circ} & 1 & \textrm{``CLUMP'' at } 24^{\circ} \\ 
\hline
2 & \textrm{``CLUMP'' at } 44^{\circ} & 2 & \textrm{``CLAMP'' at } 31^{\circ} \\ 
\hline
3 & \textrm{``RUNNY'' at } 45^{\circ} & 3 & \textrm{``CRUNK'' at } 54^{\circ} \\ 
\hline
4 & \textrm{``UNDER'' at } 49^{\circ} & 4 & \textrm{``CAMPS'' at } 64^{\circ} \\ 
\hline
5 & \textrm{``CLAMP'' at } 53^{\circ} & 5 & \textrm{``RUNNY'' at } 83^{\circ} \\ 
\hline
6 & \textrm{``CAMPS'' at } 57^{\circ} & 6 & \textrm{``UNDER'' at } 90^{\circ} \\ 
\hline
\end{array}
\]
\end{table}

The word with the smallest angle is chosen as being the most representative of the set.
When considering letter location, ``CLUMP'' is found to be most representative of the set.
However, if letter location is not considered, ``CRUNK'' is the most representative of the set.
Between the two methods the only difference in procedure was how each word was expressed as a vector.
As a result, the dominant eigenvector of the matrix and resulting angles were also different.
There is no right or wrong way of representing a word as a vector, but it is important to note these two methods can produce different results.
In testing, both methods of converting a word to a vector were compared.


\section{Testing method}

The first aspect of testing requires creating the list of possible solutions given the constraints. 
This is done by interpreting the computer's response after each round to keep only the words matching all the criteria.
For example, if after the first guess a green color is awarded to the letter ``W'' in the first location, the list of possible solutions can be narrowed down to only words which begin with the letter ``W.'' 
Similarly, if the letter ``A'' has been given a gray color, the list of possible solutions can be narrowed down further to exclude all words which contain the letter ``A.'' 
A Python program was created which took as input the computer's response after each round and output the list of possible solutions. 

The second portion of testing is selecting which word to guess next. 
This was done using the rank one approximation and latent semantic indexing method to select the word most representative of the set. 

Words to use as initial guesses needed to be chosen.
To select which words to use, the rank one approximation and latent semantic indexing method was applied to the list of all $12,947$ acceptable Wordle guesses and all $2,315$ published Wordle solutions (table 2). 
Both methods of converting a word to a vector were used to produce a list of five words to test as starting guesses: ``SOARE'', ``SLATE'', ``SORES'', ``BARES'', and ``ALERT.''

\begin{table}[ht]
\caption{Origin of words chosen as starting guesses for the model.}\label{resulttab}
\renewcommand\arraystretch{1}
\noindent\[
\begin{array}{|c|c|c|}
\hline
 \textrm{Set used}&  \textrm{Not considering letter location} & \textrm{Considering letter location} \\ 
\hline
\textrm{Acceptable guesses} & \textrm{SOARE} & \textrm{SORES},\textrm{BARES}  \\ 
\hline
\textrm{Published solutions} & \textrm{ALERT} & \textrm{SLATE} \\
\hline
\end{array}
\]
\end{table}

Considering letter location and the list of acceptable guesses, ``SORES'' is found to be most representative of the set. 
However, the double occurrence of ``S'' raises concerns about the guess's effectiveness. 
The next most representative word of the list without any double letters is ``BARES.'' 
Both ``SORES'' and ``BARES'' were tested as representatives from the set of acceptable guesses considering letter location.

Each initial guess was tested using both methods of converting a word into a vector. 
A control group with uninformed play was needed to compared the method.
The control method randomly selects a first guess, uses the available information to create a list of possible solutions, and from that list randomly selects a next guess.
It is important the baseline interprets the results and creates the same list of possible solutions as the tested methods.
However, the baseline does not make an informed decision about which word to guess next. 

\section{Results}
The 11 methods were tested with a Python program simulating both sides of the game: the player and the computer. 
To remain authentic to the conditions of the official game, the $2,315$ published Wordle solutions were used as secret words to test each method.
For each trial, the number of guesses required to solve the puzzle was recorded. 
If the method was unable to solve the secret word within six guesses, the game was counted as a loss. 
If the method correctly guessed the secret word in six or fewer guesses, the game was counted as a win. 
Of the winning games, the average number of guesses needed to solve the puzzle was determined (table 3).

\begin{table}[ht]
\caption{Results from each described method}\label{resulttab}
\renewcommand\arraystretch{1.5}
\noindent\[
\begin{array}{|c|c|c|c|c|}
\hline
\multicolumn{5}{|c|}{\textrm{Control Method}} \\
\hline
 & \multicolumn{2}{|c|}{\textrm{Avg. Guesses}} & \multicolumn{2}{|c|}{\textrm{Win \%}} \\
\hline
\textrm{RANDOM} & \multicolumn{2}{|c|}{4.59} &\multicolumn{2}{|c|}{88.2 \%} \\
\hline
\hline
\multicolumn{5}{|c|}{\textrm{Rank One Approximation with Latent Semantic Indexing}} \\
\hline
  &  \multicolumn{2}{|c|}{\textrm{Not Considering Letter Location}} & \multicolumn{2}{|c|}{\textrm{Considering Letter Location}} \\ 
\hline
 \textrm{Starting Word} & \textrm{Avg. Guesses} & \textrm{Win \%} &  \textrm{Avg. Guesses} & \textrm{Win \%} \\
\hline
 \textrm{SOARE} & 4.22 & 93.4 \% & 4.13 & 97.8 \%  \\ 
\hline
 \textrm{ALERT} & 4.21 & 96.2 \%  & 4.10 & 98.1 \%  \\ 
\hline
 \textrm{SORES} & 4.40 & 93.1 \%  & 4.26 & 97.5 \%\\ 
\hline
 \textrm{BARES} & 4.27 & 95.9 \% & 4.14 & 98.3 \% \\ 
\hline
 \textrm{SLATE} & 4.15 & 96.3 \%  & 4.04 & 98.7 \%\\ 
 \hline
\end{array}
\]
\end{table}


For every starting word, the rank one approximation and latent semantic indexing method outperformed random selection. 
The strategic edge the method provides improved the win percentage by up to $10.5\%$ and reduced the average number of guesses required by up to $0.55$.
The best strategy used ``SLATE'' as the initial guess, having both the highest win percentage and the fewest average guesses per game. 
Across all the possible starting words, no starting word dramatically outperformed the others.
Only a slight improvement of $5.6 \%$ in win percentage and $0.36$ guesses per game was seen between the best and worst starting word.
Considering the two methods of converting a word to a vector, the method which considered letter location outperformed the method which did not consider letter location for each starting word.

These results show a clear advantage to the mathematical model described in this paper over uninformed word selection.
By using a rank one approximation in conjunction with latent semantic indexing, informed guesses were made to increase the win percentage and reduce the average number of guesses to solve the puzzle.


\section{Discussion}
A Google Collab document was created with the Python code used to run these tests and calculations \cite{collab}.
Readers are invited to use the code and test alternative strategies for playing Wordle.
Other strategies could include different ways of converting a word to a vector, introducing weighted parameters to favor some letters over another, ignoring words with repeated letters, etc.

Long time Wordle players can compare their average number of guesses and win percentage to that of the method. 
For many players, their average number of guesses is better than the math model's.
How can this be if the method has perfect access to the dictionary of acceptable guesses and creates a perfect list of possible solutions?
A human player's intuition and judgement can help outperform the indiscriminate model.
Oftentimes, the program will suggest guessing a word such as ``THROE'' over ``TROPE'' if the calculations work out that the former is more representative of the set.
A human player is far more likely to correctly choose ``TROPE'' as a guess over ``THROE'' in a belief that the game chooses words which are relatively common and reasonable.

This paper presented a method of playing the game Wordle using rank one approximations and latent semantic indexing, but this general structure has applications beyond Wordle. 
Given a population with characteristics that can be represented in vector form, this method can select a member of the population most representative of the group.
As a method of surveying a group, participants which are most representative of a group could be chosen from the population.
The term ``Average Joe'' could be given merit by finding an individual who is closest to an idealized representative.
Anytime a representative from a large group is needed, these two methods provide a systematic structure of making a selection.


\end{document}